\documentclass{siamart220329}

\usepackage{amsmath}
\usepackage{amssymb}
\usepackage{amstext}
\usepackage{amscd}
\usepackage{arydshln}
\usepackage{color}
\usepackage{graphicx}
\usepackage{longtable}
\usepackage{calc}
\usepackage{microtype}
\usepackage{xcolor}
\usepackage{comment}
\usepackage{subcaption}
\usepackage{soul}

\newcommand{\diag}[1]{\mbox{${\rm diag}(#1)$}}

\newcommand{\adots}{.\:\!\raisebox{.6ex}{.}\:\!\raisebox{1.2ex}{.}}

\newcommand{\rank}[1]{\mbox{${\rm rank}(#1)$}}

\def\rank{\mathop{\rm rank}\nolimits}

\newtheorem{remark}[theorem]{{\sc Remark}}
\newtheorem{prop}[theorem]{{\sc Proposition}}
\newtheorem{example}[theorem]{{\sc Example}}

\newcounter{parenumi}
  {\end{list}}

\def\F{\mathbb{F}}

\def\R{\mathbb{R}}
\def\C{\mathbb{C}}

\def\la{\lambda}

\DeclareMathOperator{\rev}{rev}

\newlength{\tablecolwidth}
\definecolor{brilliantrose}{rgb}{1.0, 0.33, 0.64}
\definecolor{myviolet}{rgb}{0.21, 0.0, 0.85}
\definecolor{amethyst}{rgb}{0.6, 0.4, 0.8}
\definecolor{carrotorange}{rgb}{0.93, 0.57, 0.13}

\numberwithin{equation}{section}

\begin{document}

\title{Para-Hermitian rational matrices\thanks{The work of all the authors has been partially supported by the Agencia Estatal de Investigaci\'on of Spain MCIN/AEI/10.13039/501100011033/ through grants PID2019-106362GB-I00 and RED2022-134176-T.}}

\author{Froil\'an Dopico\thanks{Departamento de Matem\'aticas, Universidad Carlos III de Madrid, Avda. Universidad 30, 28911 Legan\'{e}s, Spain (\email{dopico@math.uc3m.es}).}
\and
Vanni Noferini\thanks{Department of Mathematics and Systems Analysis, Aalto University,  P.O. Box 11100, FI-00076, Finland (\email{vanni.noferini@aalto.fi , maria.quintanaponce@aalto.fi}).}
\and
Mar\'ia C. Quintana\footnotemark[3]
\and
Paul~Van~Dooren\thanks{Department of Mathematical Engineering, Universit\'{e} catholique de Louvain, Avenue Georges Lema\^itre 4, B-1348 Louvain-la-Neuve,  Belgium (\email{paul.vandooren@uclouvain.be}).}
}

\date{}

\maketitle

\begin{abstract}
In this paper we study para-Hermitian rational matrices and the associated structured rational eigenvalue problem (REP). Para-Hermitian rational matrices are square rational matrices that are Hermitian for all $z$ on the unit circle that are not poles. REPs are often solved via linearization, that is, using matrix pencils associated to the corresponding rational matrix that preserve the spectral structure. Yet, non-constant polynomial matrices cannot be para-Hermitian. Therefore, given a para-Hermitian rational matrix $R(z)$, we instead construct a $*$-palindromic linearization for $(1+z)R(z)$, whose eigenvalues that are not on the unit circle preserve the symmetries of the zeros and poles of $R(z)$. This task is achieved via M\"{o}bius transformations.
We also give a constructive method that is based on an additive decomposition into the stable and anti-stable parts of $R(z)$. Analogous results are presented for para-skew-Hermitian rational matrices, i.e., rational matrices that are skew-Hermitian upon evaluation on those points of the unit circle that are not poles.
\end{abstract}




\section{Introduction}

Rational matrices $R(t)$ play a fundamental role in systems and control theory \cite{Kai80,Rosen}, where they typically represent the transfer function of a linear time invariant system. In some important applications they have certain symmetries -- called self-conjugacy -- that need to be preserved for solving the underlying problem. For instance, such self-conjugate matrices may represent the  spectral density function of a stochastic process, or play an important role in problems such as spectral factorization, Wiener filtering or optimal control. The most common instances of self-conjugate properties that one finds in the literature are rational matrices $R(t)$ that are Hermitian (or skew-Hermitian) for every point $t$ (that is not a pole of $R(t)$) on one of the following three curves: the real axis $\R$, the imaginary axis $i\R$ and the unit circle $S^1$. To emphasize the distinction between these three cases, we will use the variable $x$ in the real axis case, the variable $s$ in the imaginary axis case and the variable $z$ in the unit circle case; when discussing results that hold for every arbitrary rational matrix, we will use $t$ as the generic variable.

In this paper, we focus on para-Hermitian and para-skew Hermitian rational matrices. To be more precise, these are square rational matrices that satisfy the following entry-wise properties (see, e.g., \cite{BN} for the para-Hermitian case):
\begin{definition} \label{def.rationals}  Let $\mathbb{F} \subseteq \mathbb{C}$ be a field. A rational matrix $R(z) \in \mathbb{F}(z)^{m\times m}$ is para-Hermitian (resp., para-skew-Hermitian) if, for all $i$ and $j$,
\begin{equation}\label{eq:giovannivanni}
    R(z)_{ij} = \overline{R(\overline{z}^{-1})_{ji}} \quad \quad (\mathrm{resp.,} \quad R(z)_{ij} = -\overline{R(\overline{z}^{-1})_{ji}} \,).
\end{equation}
If $R(z)$ satisfies the property in \eqref{eq:giovannivanni}, we shall write
\begin{equation*}R^*(z) = R\left(1/z \right)\quad \quad (\mathrm{resp.,} \quad R^*(z) = -R\left(1/z \right) \,).
\end{equation*}
\end{definition}

Para-Hermitianity is one of various possible self-conjugate properties of rational matrices. Other common such properties are mentioned in Table \ref{tab:structures}, where we also recall the name associated with each property. Note that, in the notation $R^*(t)$, the superscript $^*$ means that the coefficients of $R(t)$ are implicated in the Hermitian conjugation, but not the variable $t$; hence, $R^* (t) := [R(\overline{t})]^*$. Depending on the self-conjugate property of each class of rational matrices $R(t)$, it is easy to see that its poles and zeros will be distributed with certain symmetries, that can be described respectively as mirror images with respect to the real line $\R$, the imaginary axis $i\R$ or the unit circle $S^1$.

\begin{table}[h]
   		\begin{center}
   			\renewcommand{\arraystretch}{1.3}
   			\begin{tabular}{| c | c | c |}
   				\hline
   				Symmetry with respect to   & Hermitian: & Skew-Hermitian:  \\

         the real line $\R$ &  $R^*\left(x \right) = R(x)$ &  $R^*\left(x \right) = -R(x)$ \\ \hline
   				
   				Symmetry with respect to  & $*$-even:   &  $*$-odd:  \\

         the imaginary axis $i\R$ &  $R^*\left(s \right) = R(-s)$ &  $R^*\left(s \right) = -R(-s)$ \\ \hline

   				Symmetry with respect to  & Para-Hermitian:   & Para-skew-Hermitian:   \\

         the unit circle $S^1$ &  $R^*(z) = R\left(1/z \right)$  & $R^*(z) = -R\left(1/z \right)$  \\\hline

   			\end{tabular}

   			\caption{The structured rational matrices studied in this paper.}
   			\label{tab:structures}
   		\end{center}
   	\end{table}

In particular, the poles and zeros of a para-Hermitian (or para-skew-Hermitian) rational matrix $R(z)$ which are not on the unit circle $S^1$, appear in pairs $(\lambda,1/\overline{\lambda})$, symmetric with respect to $S^1$. In other words, each such pair have the same phase and reciprocal modulus. Moreover, if the coefficients of $R(z)$ are real then the poles and zeros that are not real, also come in complex conjugate pairs $(\lambda,\overline{\lambda})$, implying that they come in quadruples $(\lambda,1/\overline{\lambda},\overline{\lambda},1/\lambda)$ if they are not on $S^1$. This property is very important and ought to be preserved when computing the poles and zeros of $R(t)$.
Para-Hermitian matrices, not necessarily rational, that are analytic on the unit circle are extremely relevant in signal processing \cite{BN,WPP}; see in particular \cite{WPP} and the references therein for an extensive survey of their applications in various fields of engineering.

 \begin{remark} {\rm
     There is some inconsistency in the literature about the name para-Hermitian. The most prevalent custom (see, e.g., \cite{BN,WPP} and the references therein) is to call simply ``para-Hermitian" those matrices that depend on a variable $z$ and evaluate to a Hermitian matrix on the unit circle (except possibly on their singularities). However, in some papers \cite{Geninetal} these are called ``discrete-time para-Hermitian" or ``para-Hermitian on the unit circle", to distinguish them with matrices that depend on a variable $s$ and are Hermitian on the imaginary axis (except possible singularities), which are called ``continuous-time para-Hermitian" or ``para-Hermitian on the imaginary axis". In this paper, we opt for the first option which is both less pedantic and more common in the literature; we use instead the word ``even rational matrix" for the case of Hermitianity on the imaginary axis, analogously to what is done in the polynomial case \cite{AhmadMehrmann}. We point out that $*$-even rational matrices are called para-Hermitian in \cite{dopquinvan2022}; this clarification can be helpful to readers as some of the results in \cite{dopquinvan2022} will be used in some proofs.}
 \end{remark}

  Poles and zeros of rational matrices can be computed via linearizations, which ideally should preserve the structure of the original rational matrix. In Section \ref{sec:preliminaries} we review some facts about the structural data of rational matrices and of their linear system matrices. In particular, we recall the notion of strongly minimal linearization \cite{DMQVD,DQV,dopquinvan2022}. In Section \ref{sec_impossible}, we show that it is impossible to linearize a para-Hermitian rational matrix $R(z)$ with a $*$-palindromic linear system matrix. We instead show, in Section \ref{sec:Mobius}, how to construct strongly minimal $*$-palindromic linearizations for the rational matrix $(1+z)R(z)$, via M\"{o}bius transforms, if $R(z)$ is para-Hermitian. In Section \ref{sec:stable_antistable}, we introduce an additive decomposition for rational matrices in terms of their stable and anti-stable parts, and show how to construct $*$-palindromic strongly minimal linearizations for $(1+z)R(z)$ based on this decomposition. Finally, in Section \ref{sec:Taylor_and_PFD}, we pay particular attention to the construction of $*$-palindromic strongly minimal linearizations from the Taylor expansion around infinity (Subsection \ref{subsec:Taylor}) and  from the partial fraction decomposition (Subsection \ref{subsec:PFD}) of the stable part of $R(z)$. Some conclusions and lines of future research are discussed in Section \ref{sec.conclusions}.

 \section{Preliminaries}\label{sec:preliminaries}

Below, $\F$ denotes either the field of real numbers $ \mathbb{R}$ or the field of complex numbers $\mathbb{C}$; $\mathbb{F}[t]^{m\times m}$ and $\mathbb{F}(t)^{m\times m}$ denote the sets of $m \times m$ matrices whose entries are in the ring of polynomials $\mathbb{F}[t]$ and in the field of rational functions $\mathbb{F}(t)$, respectively. The elements of $\mathbb{F}[t]^{m\times m}$ and $\mathbb{F}(t)^{m\times m}$ are called, respectively, polynomial matrices (or matrix polynomials) and rational matrices. Polynomial matrices can have finite zeros (i.e., finite eigenvalues) but no finite poles, while a rational matrix $R(t)\in\mathbb{F}(t)^{m\times m}$ can have both finite poles and zeros. These are defined as follows via the local Smith--McMillan form \cite{AmMaZa15,DNZ,vandooren-laurent-1979}. Given $\la\in  \mathbb{C}$, there exist rational matrices $M_{\ell}(t)$ and $M_r(t)$ invertible at $\la$\footnote{A rational matrix $M(t)\in \F(t)^{m\times m}$ is said to be invertible at $\la\in\C$ if the constant matrix
$M(\la)$ is bounded (i.e., $M(\la)\in\C^{m\times m}$) and invertible.} such that
 \begin{equation} \label{snf}
 M_{\ell}(t)R(t)M_r(t) = \diag{(t-\la)^{d_{1}} ,\ldots, (t-\la)^{d_{r}}, 0_{(m-r)\times (m-r)}}  ,
 \end{equation}
 where $d_1\le d_2 \le \cdots \le d_r$ are integers and $r$ is the normal rank of $R(t)$. The diagonal matrix in \eqref{snf} is unique and is called the {\em local Smith-McMillan form} of $R(t)$ at $\la$ (see e.g. \cite{DMQVD}). The exponents $d_i$ are called the
 {\em structural indices or invariant orders} of $R(t)$ at $\la$. If there are strictly positive indices $0 < d_p\le \cdots \le d_r$ in \eqref{snf}, then $\la$ is a zero of $R(t)$ with partial multiplicities $(d_p, \ldots , d_r)$. If there are strictly negative indices $d_1\le \cdots \le d_q <0$ in \eqref{snf}, then $\la$ is a pole of $R(t)$ with partial multiplicities $(-d_1, \ldots , -d_q)$. The structural indices of $R(t)$ at infinity are defined as those of $R(1/t)$ at zero. The list of structural data of a rational matrix is not only formed by its finite and infinite pole and zero structures but also by its left and right minimal indices. Minimal bases and indices were introduced by Forney in \cite{For75}. The minimal bases and minimal indices of a rational matrix $R(t)$ are those associated with the following rational vector subspaces
  \[
  \begin{array}{l}
  \mathcal{N}_r (R)=\{x(t)\in\mathbb{C}(t)^{m\times 1}: R(t)x(t)=0\}, \text{and}\\
  \mathcal{N}_\ell (R)=\{y(t)^{T}\in\mathbb{C}(t)^{1\times m}: y(t)^T R(t)=0\},
  \end{array}
  \]
   which are called the right and left null-spaces of $R(t)$, respectively. If $R(t)$ is singular, then these null-spaces are non-trivial. If ${\cal N}_{r}(R)$ (resp., ${\cal N}_{\ell}(R)$) is non-trivial, it has minimal bases and minimal indices, which are called the right (resp., left) minimal bases and minimal indices of $R(t)$.

 Rosenbrock's polynomial system matrices associated with a rational matrix $R(t)$ contains its pole and zero information, whenever minimality conditions are satisfied \cite{DMQVD, DNZ, Rosen}. In order to linearize a rational matrix $R(t)$, one can thus consider linear polynomial system matrices associated with $R(t)$. These are block partitioned pencils of the form
\begin{equation}\label{eq:linearsystemmat}
L(t):= \left[\begin{array}{ccc} -t A_1 +A_0 &  t B_1 - B_0 \\ t C_1 -C_0 &  t D_1-D_0 \end{array}\right] =: \begin{bmatrix}
-A(t) & B(t)\\
C(t) & D(t)
\end{bmatrix} ,
\end{equation}
where $A(t)$ is square and regular, i.e., $\det A(t) \neq 0$, and
\begin{equation}\label{eq:transferfun}
    R(t) =D(t)+C(t)A(t)^{-1}B(t),
\end{equation}
i.e., the Schur complement of $A(t)$ in $L(t)$ is $R(t)$. Then, $R(t)$ is said to be the transfer function of $L(t).$ One can obtain the finite pole and zero structure of $R(t)$ from the eigenvalue structures of the polynomial matrices $A(t)$ and $L(t)$, respectively, provided $L(t)$ is {\em irreducible} or {\em minimal}, meaning that the matrices
\begin{equation}  \label{minimal}
\left[\begin{array}{ccc} -A(t) & B(t)  \end{array}\right] \quad \text{and}  \quad \left[\begin{array}{ccc} -A(t) \\ C(t) \end{array}\right],
\end{equation}
have, respectively, full row and column rank for all $t_0\in \mathbb{C}$ \cite{DNZ,Rosen}. In the special case that $L(t)$ is minimal and $B(t)$, $C(t)$ and $D(t)$ are constant matrices, then \eqref{eq:transferfun} is said to be a minimal generalized state-space realization of $R(t)$. Minimal polynomial system matrices in subsets of $\mathbb{C}$ and at infinity are defined and studied in \cite{DMQVD}. It was shown in \cite{DQV} that one can recover both the finite and infinite polar and zero structure of $R(t)$ from the pencils $A(t)$ and $L(t)$ provided the pencils in \eqref{minimal} have full rank also at infinity. Minimality at $\infty$ means that the matrices
 $$ \begin{bmatrix} -A_1 & B_1 \end{bmatrix}  \quad \text{and}  \quad \begin{bmatrix} -A_1\\ C_1\end{bmatrix},$$
 have, respectively, full row and column rank. If $L(t)$ is both minimal (in $\mathbb{C}$) and minimal at $\infty$ then $L(t)$ is said to be strongly minimal \cite{DMQVD,DQV} or, also, a strongly minimal linearization of $R(t)$. Moreover, in this situation, the eigenvectors, root vectors, and minimal bases of $R(t)$ can be easily recovered from those of $L(t)$, and their minimal indices are the same \cite{dopquinvan2022,NV23}.

We aim to linearize para-Hermitian rational matrices with {\em structured} linear system matrices, as in \eqref{eq:linearsystemmat}, in such a way that the chosen structure preserves the pole and zero symmetries of the original rational matrix. However, a non-constant polynomial matrix $P(z)$ and, in particular, a non-constant pencil cannot satisfy Definition \ref{def.rationals}. That is, given a para-Hermitian rational matrix $R(z)$, it is not possible to construct a non-constant para-Hermitian linear system matrix of $R(z)$.

Given a polynomial matrix $P(t)$ of degree $d$, the \emph{reversal} matrix polynomial of $P(t)$ is $\rev P(t):=t^d P\left(1/t \right).$
Recall that a matrix polynomial $P(z)$ is \emph{$*$-palindromic} (resp., \emph{$*$-anti-palindromic}) if it satisfies { \cite{GoodVibra}}
$$\rev P^*(z)=P(z) \quad \quad (\text{resp.,}\quad \rev P^*(z)=-P(z)).  $$
Then, one could try to construct a $*$-palindromic linear system matrix of $R(z)$, whose eigenvalues also have the same symmetries, i.e., they also appear in pairs $(\lambda,1/\overline{\lambda})$, symmetric with respect to $S^1$, and to apply to such a linear system matrix a structure-preserving algorithm for computing its eigenvalues \cite{KressnerQRpal,antitriangular}.
But, as we show in the following Section \ref{sec_impossible}, it is impossible to linearize a para-(skew-)Hermitian rational matrix with an (anti)-$*$-palindromic system matrix. However, if one consider the rational matrix $H(z):=(1+z)R(z)$ instead of $R(z)$, we show then how to construct a $*$-palindromic (resp., $*$-anti-palindromic) linear system matrix for $H(z)$ when $R(z)$ is para-Hermitian (resp., para-skew-Hermitian). Ultimately, the point is that, if $R(z)$ is para-Hermitian, then $H(z)$ satisfies
\begin{equation}\label{eq:one_reversal}
\rev_1H^*(z)=H(z) \quad \text{where}\quad \rev_1 H(z):=z H\left(1/z \right),
\end{equation}
and, as we show in Theorem \ref{thm:$*$-palindromic_lsm}$\rm (a)$, the transfer function of a $*$-palindromic linear system matrix satisfies \eqref{eq:one_reversal}.
Analogously, if $R(z)$ is para-skew-Hermitian, then $H(z)$ satisfies $\rev_1H^*(z)=-H(z)$.

\section{It is impossible to linearize a para-Hermitian rational matrix with a palindromic system matrix} \label{sec_impossible}

We will show in this section that it is impossible to construct a $*$-palindromic (resp., $*$-anti-palindromic) linear system matrix of a para-Hermitian (resp., para-skew-Hermitian) rational matrix $R(z)$. First, observe that Theorem \ref{thm:$*$-palindromic_lsm}$\rm (a)$ and Corollary \ref{cor:impossible} do not assume any minimality on the structured linear system matrix. See Theorem \ref{thm:$*$-palindromic_lsm}$\rm (b)$ for a result on minimality.

\begin{theorem}\label{thm:$*$-palindromic_lsm} Consider a $*$-palindromic (resp., $*$-anti-palindromic) linear system matrix $L(z)$. Then the following statements hold:
\begin{itemize}
    \item[\rm(a)] The transfer function $H(z)$ of $L(z)$ satisfies $\rev_1H^*(z)=H(z)$\\ (resp., $\rev_1 H^*(z) \allowbreak =-H(z)$).
    \item[\rm(b)] $L(z)$ is minimal if and only if it is strongly minimal.
\end{itemize}
\end{theorem}

\begin{proof} We only discuss in detail the $*$-palindromic case, since the \!\! $*$-anti-palindromic one is analogous. A $*$-palindromic linear system matrix has the form
\[ L (z) = z \begin{bmatrix}
	-A & C\\
	B^* & D
	\end{bmatrix} + \begin{bmatrix}
	-A^* & B\\
	C^* & D^*
	\end{bmatrix},\]
	where $zA+A^*$ is regular.	For $\rm (a)$, it is easy to check that $\rev_1H^*(z)=z H^*\left(1/z \right)=H(z),$ where $  H(z) = (z D + D^*) + (z B^* + C^*) (z A + A^*)^{-1} (z C+ B).$ To prove $\rm (b),$ first note that strong minimality implies minimality by definition. For the converse, assume that $L(z)=:z L_0^*+L_0$ is minimal. In particular, the first block row and the first block column of $L(0)=L_0$ have full rank. But, taking conjugate transposes, the same property holds for $L_0^*$, implying minimality at  $\infty$ and therefore strong minimality.
 \end{proof}

\begin{corollary}\label{cor:impossible} Let $R(z)\in \mathbb{F}(z)^{m\times m}$ be a nonzero rational matrix. If $R(z)$ is para-Hermitian (resp., para-skew-Hermitian) then there is no $*$-palindromic (resp., $*$-anti-palindromic) linear system matrix whose transfer function is $R(z)$.
\end{corollary}
\begin{proof} Again we only treat the $*$-palindromic case, as the $*$-anti-palindromic one is analogous. Suppose by contradiction that $R(z)$ is the transfer function of a $*$-palindromic linear system matrix. Then, by Theorem \ref{thm:$*$-palindromic_lsm}$\rm (a)$, $zR^* \left( 1/z \right)=R(z)$. But, since $R(z)$ is para-Hermitian, $R^* \left( 1/z \right)=R(z)$, and thus $zR \left( z \right)=R \left( z \right)$, which is a contradiction since $R(z) \ne 0$.
\end{proof}

While in the following sections \ref{sec:Mobius} and \ref{sec:stable_antistable} we make statements that cover both the para-Hermitian and the para-skew-Hermitian case, we give proofs only for the para-Hermitian case, as the para-skew-Hermitian one is completely analogous and left as an exercise.

\section{Para-Hermitian rational matrices and M\"{o}bius transform}\label{sec:Mobius}

In this section we consider the following M\"{o}bius transform \cite{M4Mob,Nof12} $T$ and its inverse $T^{-1}$:
\begin{equation}\label{eq_mobius1}
T:\quad x \longmapsto z=\dfrac{i-x}{i+x},\quad\quad \text{and} \quad  T^{-1}: \quad z \longmapsto x=i \dfrac{1-z}{1+z}.
\end{equation}
Note that the M\"{o}bius transformation $T$ is minus the Cayley transform. We will use the fact that $T$ maps $x \in \R$ to $T(x) \in S^1$ and, conversely, its inverse $T^{-1}$ maps $z \in S^1$ to $T^{-1}(z) \in \R$. Given a para-Hermitian rational matrix $R(z)$, we can apply the change of variable $z=T(x)$ in $R(z)$. Namely, $$R(z) \longmapsto   R(T(x))=:G(x).$$ Then, we obtain that $G(x)$ is Hermitian, i.e., $G^*(x) =G(x)$. Analogously, if $G(x)$ is a Hermitian rational matrix, then the change of variable $x = T^{-1}(z)$ maps $G(x) \longmapsto   G(T^{-1}(z))=:R(z)$ and $R(z)$ is para-Hermitian. We formalize this discussion in Lemma \ref{lem:Hermitian}.

\begin{lemma}\label{lem:Hermitian}
A rational matrix $R(z)\in \mathbb{C}(z)^{m\times m}$ is para-Hermitian (resp., para-skew-Hermitian) if and only if $G(x):=R(T(x))\in \mathbb{C}(x)^{m\times m}$ is Hermitian (resp., skew-Hermitian), where $T$ is the M\"{o}bius transformation in \eqref{eq_mobius1}. \end{lemma}

\begin{proof} Suppose that $R(z)$ is para-Hermitian, i.e., $R^*\left( 1/z \right) = R(z)$ for all $  z \in \C.$ Then, for any $x \in \C$,
	$ G^*(x) = R^*\left( \dfrac{-i-x}{-i+x}\right) = R\left( \dfrac{i-x}{i+x}\right)  = G(x) $
	and hence $G(x)$ is Hermitian. Conversely, assume now that $G(x)$ is Hermitian, i.e., $G^*\left(x \right) = G(x)$ for all $  x \in \C$. Then, for any  $z \in \C$,
	$R^*\left( 1/z \right) =G^*\left(-i \dfrac{1-1/z}{1+1/z}\right)=G\left(i \dfrac{1-z}{1+z}\right)=R(z).$
\end{proof}

Now, we can state and prove Theorem \ref{th:paraHermitian}, which is one of the main results of this paper. Before proceeding with the proof, we explain the key ideas. Observe that,
given a para-Hermitian rational matrix $R(z)$, and taking into account Lemma \ref{lem:Hermitian}, we can linearize the Hermitian rational matrix $G(x):=R(T(x))$ with a Hermitian pencil $S(x)$. For that, we can construct a Hermitian strongly minimal linear system matrix $S(x)$ of $G(x)$ as in \cite{dopquinvan2022}. We can now consider the Möbius transformation $T^{-1}$, then the rational matrix $Q(z):=S(T^{-1}(z))$ must be para-Hermitian with least common denominator $(1+z)$. Finally, if we multiply $Q(z)$ by $(1+z)$, we obtain that $(1+z)Q(z)=:L(z)$ is a $*$-palindromic (see Remark \ref{rem:$*$-palindromicpencil}) linear system matrix for $(1+z)R(z)$.

\begin{remark}\rm \label{rem:$*$-palindromicpencil} Note that if a rational matrix $Q(z)$ is para-Hermitian, that is, $Q^*\left( 1/z \right) = Q(z)$, and $L(z):=(1+z)Q(z)$ is a pencil then $L(z)$ must be $*$-palindromic. Indeed, $\rev L^*(z)=L(z)$ by \eqref{eq:one_reversal}.
\end{remark}

\begin{theorem}\label{th:paraHermitian}
Let $R(z)\in \mathbb{C}(z)^{m\times m}$ be a rational matrix. $R(z)$ is para-Hermitian (resp., para-skew-Hermitian) if and only if there exists a strongly minimal $*$-palindromic (resp., $*$-anti-palindromic) linearization of $(1+z)R(z)$.
\end{theorem}

\begin{proof}  To prove the necessity, assume first
that $R(z)$ is para-Hermitian, then we follow the next steps:
\begin{enumerate}
	\item Consider the change of variable $z=T(x)$, where $T$ is the M\"{o}bius transformation in \eqref{eq_mobius1}, and set $G(x):=R(T(x))=R\left(\dfrac{i-x}{i+x}\right).$ By Lemma \ref{lem:Hermitian}, $G(x)$ is a Hermitian rational matrix.
	\item Linearize $G(x)$ with a Hermitian strongly minimal linear system matrix $$S(x)=\begin{bmatrix}
	U_{2}- x U_{1}  & V_{2}- x V_{1}\\
	V_{2}^*- x V_{1}^* & W_{2}- x W_{1}
	\end{bmatrix}=:\left[\begin{array}{cc}
	U(x) & V(x)\\
	Y(x) & W(x)
	\end{array}\right],$$
	with $U(x)\in\C[x]^ {n\times n}$ nonsingular, $Y(x)=V^*(x)$ and $U_2$, $U_1$, $W_1$, $W_2$ Hermitian matrices. Strong minimality means that
\begin{equation}\label{eq:min_finite}
\rank\begin{bmatrix} U(x_0) \\ Y(x_0)\end{bmatrix}=\rank\begin{bmatrix} U(x_0) & V(x_0) \end{bmatrix}=n\quad \text{ for all }\quad x_0\in\C,
\end{equation}
	and that
	\begin{equation}\label{eq:min_inf}
\rank\begin{bmatrix} U_1 \\ V_1^*\end{bmatrix}=\rank\begin{bmatrix} U_1 & V_1 \end{bmatrix}=n.
	\end{equation}
	 It is always possible to construct such a Hermitian pencil $S(x)$ by \cite{dopquinvan2022}.
	\item Consider the inverse M\"{o}bius transformation $x=T^{-1}(z)$ to obtain the following para-Hermitian rational matrix $Q(z)$:
	$$Q(z):=S(T^{-1}(z))=\left[\begin{array}{cc}
	U\left(\dfrac{i-iz}{1+z} \right) & V\left(\dfrac{i-iz}{1+z} \right)\\
	Y\left(\dfrac{i-iz}{1+z} \right) & W\left(\dfrac{i-iz}{1+z} \right)
	\end{array}\right]=:\left[\begin{array}{cc}
	\widetilde U(z) & 	\widetilde V(z)\\
		\widetilde Y(z) & 	\widetilde W(z)
	\end{array}\right].$$
	Moreover, $	\widetilde W(z)- 	\widetilde {Y}(z)	\widetilde U(z) ^{-1} \widetilde V(z)=R(z)$.

	\item Finally, we set
	\begin{align*}
L(z)&:=(1+z) Q(z)=\begin{bmatrix}
	(1+z)U_{2}-i(1-z) U_{1}  & (1+z)V_{2}-i(1-z)  V_{1}\\
	(1+z)V_{2}^*-i(1-z)  V_{1}^* & (1+z)W_{2}-i(1-z) W_{1}
\end{bmatrix}
	\end{align*}

	 which is a $*$-palindromic (see Remark \ref{rem:$*$-palindromicpencil}) linear system matrix for $(z+1)R(z)$. In addition, $L(z)$ is strongly minimal. To see that $L(z)$ is strongly minimal, consider first any $\la\in\C$ with $\la\neq -1$, then by \eqref{eq:min_finite} we have that
	 	$$\rank\begin{bmatrix} 	(1+\la)U_{2}-i(1-\la) U_{1}  \\ (1+\la)V_{2}^*-i(1-\la)  V_{1}^*\end{bmatrix}=\rank\begin{bmatrix} 	U\left(i\frac{1-\la}{1+\la} \right) \\ {Y}\left(i\frac{1-\la}{1+\la} \right)\end{bmatrix} =n,$$
	 	and
	 		\begin{align*}& \rank\begin{bmatrix} 	(1+\la)U_{2}-i(1-\la) U_{1}  & { (1+\la)V_{2} -i(1-\la)  V_{1} }\end{bmatrix}\\ & =\rank\begin{bmatrix} 	U\left(i\frac{1-\la}{1+\la} \right)  & V\left(i\frac{1-\la}{1+\la} \right)\end{bmatrix} =n,	\end{align*}
	 	
	 	so that $L(z)$ is minimal at $\la$. Minimality at $-1$ follows from the fact that $S(x)$ is minimal at $\infty$, that is, by \eqref{eq:min_inf}. Finally, minimality at $\infty$ follows from the fact that $S(x)$ is minimal at $-i$.
	 	
\end{enumerate}
For the para-skew-Hermitian case, we obtain that $G(x)$ is skew-Hermitian by Lemma \ref{lem:Hermitian}. Then, we instead construct a skew-Hermitian strongly minimal linearization $\widetilde S(x) $. For that, we can use \cite{dopquinvan2022}. Finally, by following the steps above, we obtain a strongly minimal $*$-anti-palindromic linear system matrix of $(1+z)R(z).$

To prove the sufficiency, let $L(z)$ be a strongly minimal $*$-palindromic (resp., $*$-anti-palindromic) linearization of $(1+z) R(z) =: H(z)$. Then, Theorem \ref{thm:$*$-palindromic_lsm}(a) implies that $z \left(1+\dfrac{1}{z}\right) R^*\left(\dfrac{1}{z}\right) = \pm (1+z) R(z)$, where $+$ stands for the palindromic case and $-$ for the anti-palindromic. Thus,  $R^*\left(\dfrac{1}{z} \right) = \pm R(z)$.
\end{proof}

Since in Theorem \ref{th:paraHermitian} we obtain a linear system matrix $L(z)$ for $(z+1)R(z)$ instead of $R(z)$, it is important to know how to recover the  structural data of $R(z)$ from those of $L(z)$. The minimal indices are the same \cite[Theorem 2.9]{dopquinvan2022}, since multiplication by $(z+1)$ does not change the left and right rational null-spaces of $R(z)$. For the same reason the minimal bases can be recovered as described in \cite[Theorem 2.10]{dopquinvan2022}. The recovery of the
invariant orders requires some analysis. For any finite $\la\neq -1$ the invariant orders of $(z+1)R(z)$ and those of $R(z)$ at $\la$ are the same and are related to those of $L(z)$ as in \cite[Theorem 3.5]{DMQVD}. Moreover, it is easy to recover the invariant orders at $-1$ and at $\infty$ of $R(z)$ from those of $L(z)$ as we state in the following Proposition \ref{prop:invariantorders}.

\begin{prop}\label{prop:invariantorders} Let $R(t)\in\C(t)^{m\times n}$ be a rational matrix with normal rank $r$ and let $$L(t):=\left[\begin{array}{cc}
	-A(t) & B(t)\\
	C(t) & D(t)
	\end{array}\right]\in\C[t]^{(p+m)\times (p+n)},$$ with $A(t)$ nonsingular, be a linear system matrix of $(1+t)R(t)$. Then the following statements hold:
	\begin{itemize}
	\item[\rm(a)] Assume that $L(t)$ is minimal at $-1$. Let $d_1\leq\cdots\leq d_s$ be the partial multiplicities of $A(t)$ at $-1$ and let $\widetilde{d}_1\leq\cdots\leq\widetilde{d}_u$ be the partial multiplicities of $ L(t)$ at $-1.$ Then:
\begin{itemize}
		\item[\rm(a.1)] The invariant orders at $-1$ of $(1+t)R(t)$ are
		\begin{equation*}
		(-d_{s},-d_{s-1},\ldots ,-d_1,\underbrace{0,\ldots,0}_{r-s-u},\widetilde{d}_{1}, \widetilde{d}_{2},\ldots, \widetilde{d}_{u}).
		\end{equation*}
		\item[\rm(a.2)]The invariant orders of $R(t)$ at $-1$ are
		$$(-d_{s},-d_{s-1},\ldots ,-d_1,\underbrace{0,\ldots,0}_{r-s-u},\widetilde{d}_{1}, \widetilde{d}_{2},\ldots, \widetilde{d}_{u}) - (1,1,\ldots, 1).$$
\end{itemize}
	\item[\rm(b)] Assume that $L(t)$ is minimal at $\infty$. Let $e_1\leq\cdots\leq e_s$ be the partial multiplicities of $ \rev_1 A(t)$ at $0$ and let $\widetilde{e}_1\leq\cdots\leq\widetilde{e}_u$ be the partial multiplicities of $ \rev_1 L(t)$ at $0.$ Then:
	\begin{itemize}
		\item[\rm(b.1)] The invariant orders of $(1+t)R(t)$ at $\infty$ are
		\begin{equation*}
		(-e_{s},-e_{s-1},\ldots ,-e_1,\underbrace{0,\ldots,0}_{r-s-u},\widetilde{e}_{1}, \widetilde{e}_{2},\ldots, \widetilde{e}_{u}) - (1,1,\ldots, 1).
		\end{equation*}
		\item[\rm(b.2)] The invariant orders of $R(t)$ at $\infty$ are
		$$(-e_{s},-e_{s-1},\ldots ,-e_1,\underbrace{0,\ldots,0}_{r-s-u},\widetilde{e}_{1}, \widetilde{e}_{2},\ldots, \widetilde{e}_{u}).$$
	\end{itemize}
	\end{itemize}

\end{prop}

\begin{proof}
Statement $\rm(a.1)$ follows from \cite[Theorem 3.5]{DMQVD}. For $\rm(a.2)$, we consider the local Smith--McMillan form at $-1$ of $(1+t)R(t)$. That is, there exist rational matrices $M_1(t)$ and $M_2(t)$  invertible at $-1$ such that
\begin{equation} \label{eq:local-1}
(1+t)R(t) = M_1(t) \diag{(t+1)^{q_{1}} ,\ldots, (t+1)^{q_{r}}, 0_{(m-r)\times (n-r)}}M_2(t)  ,
\end{equation}
where $q_i$, for $i=1,\ldots,r$, are the invariant orders at $-1$ of $(1+t)R(t)$. If we divide \eqref{eq:local-1} by $(1+t)$, we obtain that the invariant orders at $-1$ of $R(t)$ are $q_i-1$, for $i=1,\ldots,r.$ This, together with $\rm(a.1)$, proves $\rm(a.2)$.
Statement $\rm(b.1)$ follows from \cite[Theorem 3.13]{DMQVD}. For $\rm(b.2)$, we consider \cite[Theorem 3.5]{DMQVD} applied to the linear system matrix $ \rev_1 L(t)$ that is minimal at $0$, and we obtain that
$(-e_{s},-e_{s-1},\ldots ,-e_1,\underbrace{0,\ldots,0}_{r-s-u},\widetilde{e}_{1}, \widetilde{e}_{2},\ldots, \widetilde{e}_{u})$ are the invariant orders of $t(1+1/t)R(1/t)=(1+t)R(1/t)$ at $0$. Finally, note that the invariant orders at $0$ of $R(1/t)$ are equal to those of $(1+t)R(1/t)$, since $(1+t)$ is invertible at $0$.
\end{proof}

\subsection{Alternative approaches: The choice of the M\"{o}bius transform}

We note that the use of the M\"{o}bius transform $T$ in \eqref{eq_mobius1} involves complex arithmetic when the rational matrix $R(z)$ has real coefficients. To avoid this, we can instead consider the M\"{o}bius transform $B$ and its inverse $B^{-1}$:
\begin{equation}\label{eq_bilinear}
B:\quad s \longmapsto z=\dfrac{1+s}{1-s},\quad \quad \text{and}\quad B^{-1}: \quad z \longmapsto s=\dfrac{z - 1}{z + 1}.
\end{equation}

The M\"{o}bius transformation $B$ is also called bilinear transformation (see \cite[Section 4.3.3]{Antoulas}) and maps the open left half of the complex plane onto the inside of the unit disk and the imaginary axis onto the unit circle. Given the transfer function of a discrete-time system, we can use the bilinear transformation to obtain the transfer function of a continuous-time system. In addition, if the transfer function is para-Hermitian, after the transformation it will be $*$-even as we prove in Lemma \ref{lem:para-Hermitian}.

\begin{lemma}\label{lem:para-Hermitian}
	A rational matrix $R(z)\in \mathbb{C}(z)^{m\times m}$ is para-Hermitian (resp., para-skew-Hermitian) if and only if $G(s):=R(B(s))\in \mathbb{C}(s)^{m\times m}$ is $*$-even (resp., $*$-odd), where $B$ is the M\"{o}bius transformation in \eqref{eq_bilinear}.
\end{lemma}

\begin{proof}
	If $R(z)$ is para-Hermitian, i.e., $R^*\left( 1/z \right) = R(z)$ for all $  z \in \C.$ Then, for any $s \in \C$,
	$ G^*(s) = R^*\left(  \dfrac{1+s}{1-s}\right) = R\left( \dfrac{1-s}{1+s}\right)  = G(-s) $. Conversely, if $G(s)$ is  $*$-even, i.e., $G^*\left(s \right) = G(-s)$ for all $  s \in \C$. Then, for any  $z \in \C$,
	$R^*\left( z \right) =G^*\left(\dfrac{z - 1}{z + 1}\right)=G\left( \dfrac{-z + 1}{z + 1} \right)=G\left( \dfrac{-1 + 1/z}{1 + 1/z} \right)=R(1/z).$
\end{proof}

\begin{remark} \rm
    Theorem \ref{th:paraHermitian} can be also proved by using the bilinear transformation $B$ instead of minus the Cayley transform $T$, since $*$-even (resp., $*$-odd ) rational matrices admit $*$-even (resp., $*$-odd ) linearizations \cite{dopquinvan2022}.
\end{remark}

 \begin{remark} \rm Both proofs of Theorem \ref{th:paraHermitian} are constructive so that one can obtain a strongly minimal $*$-palindromic (resp., $*$-anti-palindromic) linear system matrix of $(1+z)R(z)$ from it. However, the linear system matrices that one obtains from the proofs are not unique and will depend on the linearization that one constructs after the corresponding Möbius transformation for the Hermitian (resp., skew-Hermitian) rational matrix $G(x)$ or for the $*$-even (resp., $*$-odd ) rational matrix $G(s)$. As we mentioned above, we can consider the linearizations for Hermitian (resp., skew-Hermitian) or $*$-even (resp., $*$-odd ) rational matrices constructed in \cite{dopquinvan2022}, but this implies computing the matrix coefficients of the Laurent expansion around infinity of $G(x)$ or $G(s)$. In the next sections, we give alternative methods that do not involve such computations when the corresponding rational matrix has no poles on the unit circle.
\end{remark}

To optimize the accuracy of numerical algorithms, it might also be desirable to make a choice of the M\"{o}bius transform so that the invariant orders at $\la=-1$ are not ``shifted", as it is explained in Proposition \ref{prop:invariantorders}, if $\la=-1$ is a pole or a zero of $R(z)$. To avoid this, we can instead consider the M\"{o}bius transform $B_{\alpha}$ and its inverse $B_{\alpha}^{-1}$:
\begin{equation}\label{eq_transf}
B_{\alpha}:\quad s \longmapsto z=\dfrac{\alpha+\alpha s}{\overline{\alpha}-\overline{\alpha} s}, \quad \quad \text{and}\quad B_{\alpha}^{-1}: \quad z \longmapsto s=\dfrac{\overline{\alpha} z - \alpha}{\overline{\alpha} z + \alpha},
\end{equation}
with $\alpha\in\C$,  $\alpha \ne 0$. Then, analogous to Lemma \ref{lem:para-Hermitian}, we have the following Lemma \ref{lem:para-Hermitian2}.
\begin{lemma}\label{lem:para-Hermitian2}
    A rational matrix $R(z)\in \mathbb{C}(z)^{m\times m}$ is para-Hermitian (resp., para-skew-Hermitian) if and only if $G(s):=R(B_{\alpha}(s))\in { \mathbb{C}(s)^{m\times m}}$ is $*$-even (resp., $*$-odd), where $B_{\alpha}$ is the M\"{o}bius transformation in \eqref{eq_transf}.
\end{lemma}

 By using the M\"{o}bius transformation $B_{\alpha}$ in \eqref{eq_transf} and Lemma \ref{lem:para-Hermitian2}, we can construct a $*$-palindromic linearization for $H_{\alpha}(z):=(\alpha +\overline{\alpha} z )R(z)$, as we state in Theorem \ref{th:paraHermitian2}. Then, we can choose $\alpha$ in $B_{\alpha}$ such that $-\alpha^2/|\alpha|^2$ is neither a pole nor a zero of $R(z)$, to avoid the invariant orders of finite poles and/or finite zeros of $R(z)$ to be ``shifted". Theorem \ref{th:paraHermitian2} can be proved by following  similar steps to those in the proof of Theorem \ref{th:paraHermitian}. Note that if $R(z)$ is para-Hermitian, then $H_{\alpha}(z)$ satisfies \eqref{eq:one_reversal}, i.e., $\rev_1H_{\alpha}^*(z)=H_{\alpha}(z).$

\begin{theorem}\label{th:paraHermitian2}
Let $R(z)\in \mathbb{C}(z)^{m\times m}$ be a rational matrix and $\alpha \in \mathbb{C}, \; \alpha \ne 0$. $R(z)$ is para-Hermitian (resp., para-skew-Hermitian) if and only if there exists a strongly minimal $*$-palindromic (resp., $*$-anti-palindromic) linearization of $(\alpha +\overline{\alpha} z )R(z)$.
\end{theorem}

\subsection{Parametrizing para-Hermitian rational matrices}
The definition of para-Hermitian (resp., para-skew-Hermitian) rational matrices does not provide an explicit method for constructing all the matrices in this class. As a direct corollary of the results in this section we present in Corollary \ref{cor.parameters} different ways to generate all the matrices in this class. Observe that this result does not assume any minimality of the involved matrices.

\begin{corollary} \label{cor.parameters} Let $0\neq \alpha \in \mathbb{C}$. Then $R(z) \in \mathbb{C} (z)^{m\times m}$ is
para-Hermitian (resp., para-skew-Hermitian) if and only if there exist constant matrices $A\in \mathbb{C}^{n\times n}$, $C, B \in \mathbb{C}^{n\times m}$, and $D\in \mathbb{C}^{m\times m}$ such that the pencil $z A + A^*$ (resp., $z A - A^*$) is regular and
\begin{equation} \label{eq.parameter}
R(z) = \frac{1}{\alpha + \overline{\alpha} z} \left[ (z D + D^*) + (z B^* + C^*) (z A + A^*)^{-1} (z C+ B) \right]
\end{equation}
\begin{equation} \label{eq.2parameter} (resp., \quad
R(z) = \frac{1}{\alpha + \overline{\alpha} z} \left[ (z D - D^*) + (z B^* - C^*) (z A - A^*)^{-1} (z C- B) \right] \quad).
\end{equation}
\end{corollary}
\begin{proof}
It is easy to check that \eqref{eq.parameter} (resp., \eqref{eq.2parameter}) is para-Hermitian (resp., para-skew-Hermitian). The fact that any para-Hermitian (resp., para-skew-Hermitian) rational matrix can be written as in \eqref{eq.parameter} (resp., \eqref{eq.2parameter}) follows from Theorem \ref{th:paraHermitian2} and the structure of any $*$-palindromic (resp., $*$-anti-palindromic) pencil.
\end{proof}

\section{Decomposition into stable and anti-stable parts}\label{sec:stable_antistable}

Every rational matrix $R(t)$ can be written with an additive decomposition as in Lemma \ref{lem_decomp}, where the rational matrices $R_{in}(t)$ and $R_{out}(t)$ are called the stable and anti-stable parts of $R(t)$, respectively; here, a rational matrix is called stable (resp. anti-stable) if all its poles have moduli strictly smaller (resp. strictly larger) than $1$. If $t=z$ and $R(z)$ is para-Hermitian, the proper rational matrix $R_{p}(z)$ defined in Lemma \ref{lem_decomp} is also para-Hermitian. Then we will see  in Theorem \ref{th:poles_unitcircle} that we only need to perform a M\"{o}bius transform on $R_{p}(z)$ to construct a strongly minimal $*$-palindromic linearization $L(z)$ of $(1+z)R(z)$. In particular, if $R(z)$ has no poles on the unit circle, we can construct $L(z)$ without considering a M\"{o}bius transform as we state in Theorem \ref{th:nopoles_unitcircle}.

\begin{lemma}\label{lem_decomp} Let $R(t) \in \mathbb{C} (t)^{m \times n}$ be a rational matrix. Then there exists an additive decomposition of the form
\begin{equation}\label{eq_decomp}
R(t)= R_{in}(t) + R_{out}(t) +   R_{S^1}(t) + R_0,
\end{equation}
 where $R_{in}(t)$ is a strictly proper rational matrix that has all its poles inside the open unit disk, and is therefore stable; $R_{out}(t)$ is such that $R_{out}(0)=0$ and has all its poles (infinity included) strictly outside the unit circle, and is therefore anti-stable; $R_{S^1}(t)$ is a strictly proper rational matrix that has all its poles on the unit circle; and $R_0$ is a constant matrix. Moreover, the decomposition in \eqref{eq_decomp} is unique. In addition,  $R(z)   \in \mathbb{C} (z)^{m \times m}$ is para-Hermitian (resp., para-skew-Hermitian)  if and only if
\begin{equation}\label{eq_selfconjugate1}
 R_{in}^*(z) = R_{out}(1/z) \quad\quad  (\text{resp.,}\quad R_{in}^*(z) = -R_{out}(1/z)),
\end{equation}
and the proper rational matrix $R_p(z):=R_{S^1}(z)+R_0$ is para-Hermitian (resp., para-skew-Hermitian).

\end{lemma}

\begin{proof}  Consider the (unique) partial fraction expansion
\[ R(t) = P(t) + \sum_{\la \in \C} \sum_{j=1}^{d(\la)} \frac{1}{(t-\la)^j}A_{\la,j}, \]
where $P(t)$ is a polynomial matrix, $d(\la)$ denotes the degree of the pole $\la$, and $A_{\lambda , j} \in \mathbb{C}^{m\times n}$ are constant matrices. We can now define
\[
    R_{in}(t) := \sum_{|\la|<1} \sum_{j=1}^{d(\la)} \frac{1}{(t-\la)^j}A_{\la,j}\quad  \text{ and }\quad   R_{S^1}(t) := \sum_{|\la|=1} \sum_{j=1}^{d(\la)} \frac{1}{(t-\la)^j}A_{\la,j}. \]
    It is clear that $R_{in}(t)$ and $R_{S^1}(t)$ are strictly proper, and that the poles of $U(t):=R_{in}(t)+R_{S^1}(t)$ are those poles of $R(t)$ that are in the closed unit disk; the uniqueness of $R_{in}(t)$ and of $R_{S^1}(t)$ follows from the uniqueness of the partial fraction expansion. Moreover, $S(t):=R(t)-U(t)$ is analytic at $0$, because the coefficients of the Laurent expansion around $t=0$ with negative indices of $R(t)$ and $U(t)$ are equal by construction. Hence, we can (uniquely) define $R_0:=S(0)$ and $R_{out}(t):=S(t)-R_0$.

    We assume now that $R(z)$ is para-Hermitian. Then, since $R^*(z)=R(1/z)$, we have that
    \begin{equation}\label{eq_paraHerm_decomp}
        R_{in}^*(z) +  R_{out}^*(z) + R_p^*(z)=R_{in}(1/z)  + R_{out}(1/z) + R_p(1/z),
    \end{equation}
    where $R_p(z):=R_{S^1}(z)+R_0.$ We note that both $R_{in}^*(z)$ and $R_{out}(1/z)$ are strictly proper rational matrices that have all their poles in the open unit disk. Since the decomposition in \eqref{eq_paraHerm_decomp} for $R^*(z)$ is unique, this implies that $R_{in}^*(z)=R_{out}(1/z)$. Then, the rational matrix $H(z):=R_{in}(z)+R_{out}(z)$ is para-Hermitian and, therefore, $R_p(z)=R(z)-H(z)$ is also para-Hermitian.

Finally, let us assume that \eqref{eq_selfconjugate1} holds and that $R_p(z)$ is para-Hermitian (resp., para-skew-Hermitian). Then,
\[
R^* (z)  = R_{in}^* (z) +  R_{out}^* (z) + R_p^* (z)
         = \pm R_{out} (1/z) \pm R_{in} (1/z) \pm R_p (1/z)
         = \pm R(1/z),
\]
where $+$ stands for the para-Hermitian case and $-$ for
the para-skew-Hermitian. This concludes the proof.
\end{proof}

Next, we assume that $R(z)$ is para-Hermitian (resp., para-skew-Hermitian) and has no poles on the unit circle. That is, the additive decomposition in \eqref{eq_decomp} is of the form
\begin{equation}\label{eq_decomp2}
R(z)= R_{in}(z)  + R_{out}(z) +  R_0, \quad\text{ with }\, R_{in}^*(z)=R_{out}(1/z) \,\text{ and }\, R_0^*=R_0,
\end{equation}
$$ (resp., \;
R(z)= R_{in}(z)  + R_{out}(z) +  R_0, \quad\text{ with }\, R_{in}^*(z)= - R_{out}(1/z) \,\text{ and }\, R_0^*=-R_0 ).
$$
Then, we have the following result.
\begin{theorem}\label{th:nopoles_unitcircle}
	Let $R(z)\in \mathbb{C}(z)^{m\times m}$ be a para-Hermitian (resp., para-skew-Hermitian) rational matrix having no poles on the unit circle. Consider an additive decomposition of $R(z)$ as in \eqref{eq_decomp2} such that $R_{S^1}(z)=0$, and a minimal generalized state-space realization of $R_{in}(z)$: \begin{equation} \label{eq:Rin1strep}
R_{in}(z) = B(  z A_1 - A_0)^{-1}C,
 \end{equation}
	with $A_1$ invertible. Then,
	$$R_{out}(z)= z C^*( A_1^* -  z A_0^* )^{-1}B^*,\quad\quad(\text{resp.,}\quad R_{out}(z)= z C^*(  z A_0^*-A_1^*  )^{-1}B^*)$$
	is a minimal generalized state-space realization of $R_{out}(z)$, and the following pencil $L(z)$ is a strongly minimal $*$-palindromic (resp., $*$-anti-palindromic) linearization of $(1+z)R(z)$:
	$$   L(z) = \left[\begin{array}{cc|c} 0 & A_0 -  z A_1  & C \\ z A_0^* -  A_1^*   & 0 & B^*(1+z) \\ \hline z C^* & B(1+z) & R_0(1+z) \end{array}\right] .
	$$
 $$(\text{resp.,}\quad L(z) = \left[\begin{array}{cc|c} 0 & A_0 -  z A_1  & C \\    A_1^* - z A_0^*  & 0 & -B^*(1+z) \\ \hline -z C^* & B(1+z) & R_0(1+z) \end{array}\right] ).$$
	
\end{theorem}

\begin{proof}
 The expression of the realization of $R_{out} (z)$ follows directly from $R^*_{in} (z) = R_{out} (1/z)$ in \eqref{eq_decomp2}. The minimality of  the realization of $R_{out} (z)$ at $z\ne 0$ follows from the minimality of the realization of $R_{in} (z)$ and at $z=0$ from the invertibility of $A_1$.
 It is clear that $ L(z)$ is $*$-palindromic  and that its transfer function is $(1+z) R(z)$. Then,  the minimality of $L(z)$ at any finite $z \ne -1$ follows from the fact that $R_{in}(z)$ and $ R_{out}(z)$ have no poles in common. The minimality at $z=-1$, follows from the fact that $-1$ is not a pole of $R(z)$ by assumption. Finally, the strong minimality of $L(z)$ follows from Theorem \ref{thm:$*$-palindromic_lsm}(b).
 \end{proof}

In Section \ref{sec:Taylor_and_PFD}, we show how to construct  a strongly minimal $*$-palindromic linearization of $(1+z) R(z)$ for a para-Hermitian rational matrix $R(z)$ without poles on the unit circle from the Taylor expansion of $R_{in}(z)$ and from the partial fraction decomposition of $R_{in}(z)$, instead of directly from a minimal generalized state-space realization \eqref{eq:Rin1strep} as we did in Theorem \ref{th:nopoles_unitcircle}.

If there are poles on the unit circle, we can consider a M\"{o}bius transform for the para-Hermitian proper rational matrix $R_{p}(z)$, as we described in Section \ref{sec:Mobius}, to construct a strongly minimal $*$-palindromic linearization $L_p(z)$ for $(1+z)R_p(z)$. Then, as we  state in Theorem \ref{th:poles_unitcircle}, we can combine $L_p(z)$ with a $*$-palindromic linearization constructed from the stable and anti-stable parts, to construct a strongly minimal $*$-palindromic linearization for $(1+z)R(z)$. In Theorem \ref{th:poles_unitcircle}, whose proof is immediate and left as an exercise to the reader, we only consider for brevity the para-Hermitian case but the para-skew-Hermitian case is analogous.
\begin{theorem}\label{th:poles_unitcircle}
	Let $R(z)\in \mathbb{C}(z)^{m\times m}$ be a para-Hermitian rational  matrix expressed by the additive decomposition in Lemma \ref{lem_decomp}. Consider a strongly minimal $*$-palindromic linearization of $(1+z)[R_{in}(z)+R_{out}(z)]$:
	$$   L_{in/out}(z) := \left[\begin{array}{cc|c} 0 & A_0 -  z A_1  & C \\ z A_0^* -  A_1^*   & 0 & B^*(1+z) \\ \hline z C^* & B(1+z) & 0 \end{array}\right],
	$$
 as in Theorem \ref{th:nopoles_unitcircle}, and a strongly minimal $*$-palindromic linearization of $(1+z)R_{p}(z)$:
$$ L_{p}(z):=\left[\begin{array}{cc}
	-A_p(z) & B_p(z)\\
	C_p(z) & D_p(z)
	\end{array}\right].$$
 Then the following pencil is a strongly minimal $*$-palindromic linearization of $(1+z)R(z)$:

 \begin{equation}\label{eq:NONAME}
			L(\la):=  \left[\begin{array}{ccc|c}
			    0 & A_0 -  z A_1  & 0 & C \\
   z A_0^* -  A_1^*   & 0 & 0 & B^*(1+z)\\
       0 & 0 & -A_{p}(z) & B_{p}(z) \\ \hline \phantom{\Big|}
       z C^* & B(1+z)  & C_p(z) &    D_p(z) \end{array}\right].
\end{equation}

\end{theorem}

\section{Construction of linearizations from different representations of $R_{in}(z)$}\label{sec:Taylor_and_PFD}

In this section we assume that $R(z)$ is a para-Hermitian rational matrix having no poles on the unit circle. That is, $R_p(z)$ in the additive decomposition \eqref{eq_decomp} is a constant matrix $R_0$ and $R(z)$ can be written as in \eqref{eq_decomp2}. Otherwise, we can consider a M\"{o}bius transform for the para-Hermitian proper rational matrix $R_{p}(z)$ and use Theorem \ref{th:poles_unitcircle}.  Our first goal is to complement the construction in Theorem \ref{th:nopoles_unitcircle} by showing how to construct explicitly a minimal generalized state-space realization of $R_{in} (z)$ as in \eqref{eq:Rin1strep} from a Taylor expansion of $R_{in} (z)$ around the point at infinity. Observe in this context that there are infinitely many realizations as in \eqref{eq:Rin1strep} of any $R_{in} (z)$. In the second place, we explore explicit constructions of strongly minimal $*$-palindromic linearizations of $(1 + z)R(z)$ from the partial fraction decomposition of $R_{in} (z)$. We omit in this section to state results for para-skew-Hermitian rational matrices for brevity. They can be easily deduced as in previous sections.

\subsection{ Construction from the Taylor expansion of $R_{in}(z)$}\label{subsec:Taylor}

The strictly proper rational matrix $R_{in}(z)$ in \eqref{eq_decomp} can be represented via its Taylor expansion around the point at infinity. Namely,
\begin{equation}\label{eq:sproper}
R_{in}(z):= R_{-1} z ^{-1} + R_{-2} z ^{-2} +  R_{-3} z ^{-3} + \cdots .
\end{equation}
In this section, we construct a strongly minimal $*$-palindromic linear system matrix { as} in Theorem \ref{th:nopoles_unitcircle}, when $R_{in}(z)$ is represented as in \eqref{eq:sproper} by using the algorithm in \cite[Section 3.4]{real} or \cite[Theorem 5.1]{dopquinvan2022}. For that, let us consider the following block Hankel matrix $H$ and shifted block Hankel matrix  $H_\sigma$ associated with $R_{in}(z)$:
\begin{equation} \label{Hankel} H :=  \left[\begin{array}{cccc} R_{-1} & R_{-2}  & \ldots & R_{-k}  \\ [+2mm] R_{-2}  &  & \adots &  R_{-k-1} \\   [+2mm]
\vdots & \adots & \adots & \vdots \\  [+2mm] R_{-k} &  R_{-k-1}  & \ldots & R_{-2k+1}
\end{array}\right], \,\,  H_\sigma :=  \left[\begin{array}{cccc} R_{-2} & R_{-3}  & \ldots & R_{-k-1}  \\ [+2mm] R_{-3}  &  & \adots &  R_{-k-2} \\   [+2mm]
\vdots & \adots & \adots & \vdots \\  [+2mm] R_{-k-1} &  R_{-k-2}  & \ldots & R_{-2k}
\end{array}\right].
\end{equation}
Then for sufficiently large $k$ the rank $r_f$ of $H$ equals the total polar degree of $R_{in}(z)$, i.e., the sum of the degrees of the denominators in the Smith-McMillan form of $R_{in}(z)$ \cite{Kai80}. We assume in the sequel that we are taking such a sufficiently large $k$. The algorithm in \cite[Section 3.4]{real} for strictly proper rational matrices or \cite[Theorem 5.1]{dopquinvan2022} implies that the linear system matrix in Lemma \ref{th:lin_sproper} is a strongly minimal linearization for $R_{in}(z)$.

\begin{lemma} \label{th:lin_sproper} { {\rm \cite[Theorem 5.1]{dopquinvan2022}}}
	Let $R_{in}(z)\in \mathbb{C}(z)^{m\times n}$ be a strictly proper rational matrix as in \eqref{eq:sproper}. Let $H$ and $H_\sigma$  be the block Hankel matrices in \eqref{Hankel} and $r_f:=\rank H$. Let $U:=\left[\begin{array}{cc}U_1 & U_2\end{array}\right]$ and $V:=\left[\begin{array}{cc}V_1 &  V_2\end{array}\right]$ be unitary matrices  such that
	\begin{equation} \label{rankrf}
	U^* HV =\left[\begin{array}{cc} \widehat H &  0 \\ 0 & 0 \end{array}\right] =\left[\begin{array}{cc} U_1^* HV_1 &  0 \\ 0 & 0 \end{array}\right],
	\end{equation}
	where $\widehat H$ is $r_f \times r_f$ and invertible. Let us now partition the matrices $U_1$ and $V_1$ as follows
	\begin{equation} \label{partition} U_1=\left[\begin{array}{cc} U_{11} \\ U_{21} \end{array}\right] \quad \mathrm{and} \quad V_1=\left[\begin{array}{cc} V_{11} \\ V_{21} \end{array}\right],
	\end{equation}
	where the matrices $U_{11}$ and $V_{11}$ have dimension $m\times r_f$ and $n\times r_f$, respectively. Then
	\begin{equation} \label{realize} L_{sp}(z):= \left[\begin{array}{c|c} U_1^* H_\sigma V_1 - z \widehat H  &  \widehat H V_{11}^*  \\ \hline \phantom{\Big|} U_{11}\widehat H & 0 \end{array}\right]
	\end{equation} is a strongly minimal linearization for $R_{in}(z)$. In particular,
	\begin{equation}\label{eq_in}
R_{in}(z) =  U_{11}\widehat H(z \widehat H- U_1^* H_\sigma V_1)^{-1}\widehat H V_{11}^*.
	\end{equation}
	 \end{lemma}
 Then, we have the following result for the construction of { a} strongly minimal $*$-palindromic linear system matrix of $(1+z)R(z)$, that follows from Lemma \ref{th:lin_sproper} and Theorem \ref{th:nopoles_unitcircle} { just by replacing \eqref{eq:Rin1strep} by \eqref{eq_in}}.

 \begin{theorem}\label{th:Lurent} 	Let $R(z)\in \mathbb{C}(z)^{m\times m}$ be a para-Hermitian rational matrix having no poles on the unit circle. Consider an additive decomposition of $R(z)$ as in \eqref{eq_decomp2} and a minimal { generalized} state-space realization of $R_{in}(z)$ as in Lemma \ref{th:lin_sproper}. Namely, \begin{equation*}
 	R_{in}(z) =  U_{11}\widehat H(z \widehat H- U_1^* H_\sigma V_1)^{-1}\widehat H V_{11}^* .
 	\end{equation*}
 	Then,
 	\begin{equation}\label{eq_out}
 	R_{out}(z) = z V_{11} \widehat H^*( \widehat H^* - z V_1^*H_\sigma^* U_1)^{-1}\widehat H^* U_{11}^*
 	\end{equation}
 	is a minimal { generalized} state-space realization of $R_{out}(z)$,	and the following pencil $L(z)$ is a strongly minimal $*$-palindromic linear system matrix of $(1+z)R(z)$:
 	 $$   L(z) = \left[\begin{array}{cc|c} 0 &  U_1^* H_\sigma V_1 - z \widehat H& \widehat H V_{11}^* \\  z V_1^* H_\sigma^* U_1 - \widehat H^*  & 0 & \widehat H^* U_{11}^* (1+z) \\ \hline z V_{11} \widehat H^* & U_{11}\widehat H(1+z) & R_0(1+z) \end{array}\right].
 	$$
 \end{theorem}

\subsection{Construction from the partial fraction decomposition of $R_{in}(z)$}\label{subsec:PFD}
	
	 Let us assume that $\Omega_{in}$ is the set of poles of  the para-Hermitian rational matrix $R(z) \in \mathbb{C}(z)^{m\times m}$ inside the unit circle, that $\Omega_{out}$ is the set of poles of $R(z)$ outside the unit circle, and that $R(z)$ has no poles on the unit circle. We { also assume that} the strictly proper rational matrix $R_{in}(z)$ in \eqref{eq_decomp} { is expressed} via its partial fraction decomposition. Namely,
	\begin{equation}\label{eq:sproper_partial}
	R_{in}(z) = \sum_{\la_i \in \Omega_{in}} R_{\la_i}(z) \quad \text{ with }\quad R_{\la_i}(z):= \sum_{j=1}^{d_i} \frac{R_j}{(z-\la_i)^j},
	\end{equation}
 { where $d_i$ denotes the degree of the pole $\lambda_i$.}
	Then, since { $ R_{out}(z) = R_{in}^*(1/z)$}, we have that
		\begin{equation}\label{eq:out_partial}
	R_{out}(z) = \sum_{\la_i \in \Omega_{in}} R_{\la_i}^*(1/z) \quad \text{ with }\quad  R_{\la_i}^*(1/z)= \sum_{j=1}^{d_i} \frac{z^j R_j^*}{(1-z\overline{\la_i})^j}.
	\end{equation}

Let us first assume that $R(z)$ has only one pole $\la$ inside $S^1$. That is, { taking into account \eqref{eq_decomp2}}, $R(z)$ can be written as
\begin{equation}\label{eq: rat_one_pole}
    R(z)= R_{\la}(z) + R_0+ R_{\la}^*(1/z)\quad \text{ with }\quad R_{\la}(z)= \sum_{j=1}^{d} \frac{R_j}{(z-\la)^j}\quad {\text{ and }\quad R_0^* = R_0} .
\end{equation}
 We set:

 $$K_{\la}(z):=\begin{bmatrix}
 (\la-z) I_m &  & & &  \\
 I_m &  (\la-z) I_m & & &  \\
 &  \ddots & \ddots & & \\
 & &I_m & (\la-z) I_m \end{bmatrix}\in { \mathbb{C}[z]^{d m\times d m} },\; F:=\begin{bmatrix}
     R_{d}   \\
     R_{d-1}   \\
      \vdots \\
      R_{1}
\end{bmatrix}.$$
Then, by assuming that $R_{d}$ is invertible, we obtain the following result.

 \begin{theorem}\label{th:$*$-palindromic_onepole} Let $R(z)\in {\mathbb{C}(z)^{m\times m}}$ be a para-Hermitian rational matrix as in \eqref{eq: rat_one_pole} with { $\lambda$ inside the unit circle $S^1$ and} $R_{d}$ invertible. Then the following pencil $L_{\la}(z)$ is a strongly minimal $*$-palindromic linearization for $(1+z)R(z)$:

\begin{equation}\label{pal_pencil_anydegree}
L_{\la}(z) := \left[\begin{array}{c;{2pt/2pt}c|c}
 0 & K_{\la}(z)  & F \\ \hdashline[2pt/2pt] \phantom{\Big|}
 \rev K_{\la}^*(z)   & 0 & (1+z)(e_{d}\otimes I_m) \\ \hline \phantom{\Big|}
 z F^* & (1+z)(e_{d}^T\otimes I_m) & (1+z)R_0 \end{array}\right],
 \end{equation}
where $e_{d}$ is the last canonical vector of { $\C^{d}$} and $\otimes$ denotes the Kronecker product.
\end{theorem}

\begin{proof}  It is easy to see that $L_{\la}(z)$ is $*$-palindromic, and the strong minimality follows from the facts that $|\lambda|<1$ and $R_{d}$ is invertible. To show that the transfer function of $L_{\la}(z)$ is $(1+z)R(z)$ we can use \cite[Lemma 3.1]{dopquinvan2022}. First, we set $L_{\la}(z) = \left[\begin{array}{c|c} -A(z) & B(z) \\ \hline  C(z) & D(z) \end{array}\right].$ Now, we consider the following rational matrix:
 $$H(z):=\left[\begin{array}{ccc;{2pt/2pt}ccc|c}
\dfrac{(1+z)z^{{d}-1}}{(1-z\overline{\la})^{d}}I_m & \cdots & \dfrac{1+z}{1-z\overline{\la}}I_m & Z_{d}^T(z) & \cdots & Z_{1}^T(z)  & I_m
\end{array}\right],$$
where $Z_k(z):=\displaystyle \sum_{j=k}^{d} \frac{R_{j}}{(z-\la)^{j-k+1}}$, for $k=1,\ldots,d$. Then, by \cite[Lemma 3.1]{dopquinvan2022}, $$\left[\begin{array}{cc}C(z) & D(z)\end{array}\right]H(z)^{T}=(1+z)R(z)$$ is the transfer function of $L_{\la}(z)$, since $\left[\begin{array}{cc}-A(z) & B(z)\end{array}\right] H(z)^T=0$, { $H(z)$ has full row normal rank and the right-most block of $H(z)$ is $I_m$}. \end{proof}

We now construct an equivalent linear system matrix to $L_{\la}(z)$, that is also a strongly minimal $*$-palindromic linearization for $(1+z)R(z)$, { where $R(z)$ is as in \eqref{eq: rat_one_pole},} when $R_{d}$ is invertible. For that, we consider the block Hankel matrix:
\begin{equation} \label{eq_Hankelinout}
H:= \left[\begin{array}{cccc} & & & R_{d} \\  & & \adots &  R_{d-1} \\  & R_{d} & \adots & \vdots \\
R_{d} &  R_{d-1} & \ldots & R_{1}
\end{array}\right] { \in \mathbb{C}^{dm \times dm}},
\end{equation}
and we set
$$M_{\la}(z):=\begin{bmatrix}
   &  & &  S_{d}(z)  \\ 	
   &   & \adots & R_{d} + S_{d-1}(z) \\	
  & S_{d}(z) &  \adots &  \vdots & \\	
  S_{d}(z) & R_{d} + S_{d-1}(z) & \cdots & R_{2} + S_{1}(z)
\end{bmatrix}, \quad\quad G:=\begin{bmatrix}
     R_{d}^*   \\
     R_{d-1}^*   \\
      \vdots \\
      R_{1}^*
\end{bmatrix},$$
with $S_j(z):=(\la-z)R_j.$ Then, we obtain the following result.

\begin{theorem}\label{th_linear_palind_Hankel} Let $R(z)\in { \mathbb{C}(z)^{m\times m}}$ be a para-Hermitian rational matrix as in \eqref{eq: rat_one_pole} with { $\lambda$ inside the unit circle $S^1$ and} $R_{d}$ invertible. Then the following pencil $\widehat L_{\la}(z)$ is a strongly minimal $*$-palindromic linearization for $(1+z)R(z)$:

\begin{equation}\label{eq_linear_palind}
\widehat L_{\la}(z) = \left[\begin{array}{c;{2pt/2pt}c|c} 0 &  M_{\la}(z) &  F \\ \hdashline[2pt/2pt] \phantom{\Big|} \rev M_{\la}^*(z) & 0 & (1+z) G \\ \hline \phantom{\Big|}
z F^* & (1+z) G^*  & (1+z) R_0  \end{array}\right].
\end{equation}
\end{theorem}

\begin{proof} { Observe that $ \widehat L_{\la}(z)$ satisfies the following identity:}
    \begin{equation}\label{eq_equiv_transfer}
    \widehat L_{\la}(z)  = \left[\begin{array}{cc|c} { I_{dm}} & &  \\ & H^* &  \\ \hline  &  &  I_m \end{array}\right] L_{\la}(z) \left[\begin{array}{cc|c} {  I_{dm}} & &  \\ & H &  \\ \hline  &  & I_m \end{array}\right],
\end{equation}
where $L_{\la}(z)$ is the strongly minimal $*$-palindromic linearization in Theorem \ref{th:$*$-palindromic_onepole}.  Note that $\widehat L_{\la}(z)$ is also $*$-palindromic. In addition, since $R_{d}$ is invertible, $\widehat L_{\la}(z)$ is also a strongly minimal linearization for $(1+z)R(z)$, since the transformation in \eqref{eq_equiv_transfer} preserves the transfer function of $L_{\la}(z)$ \cite[Ch. 2, Thm. 3.1]{Rosen} and the { (strong)} minimality in that case.\end{proof}

If $R_{d}$ is not invertible, we can construct a trimmed linear system matrix from $ \widehat L_{\la}(z)$ that is a strongly minimal $*$-palindromic linearization for $(1+z)R(z)$, { as we show in Theorem \ref{th:$*$-palindromic}. The proof of Theorem \ref{th:$*$-palindromic} is inspired by the proof of \cite[Theorem 3.2]{dopquinvan2022}, a result which is valid only for (unstructured) polynomial matrices.}

\begin{theorem}\label{th:$*$-palindromic}
	Let $R(z)\in { \mathbb{C}(z)^{m\times m}}$ be a para-Hermitian rational matrix as in \eqref{eq: rat_one_pole} { with $\lambda$ inside the unit circle $S^1$ and $R_d \ne 0$}. Let $H$ be the block Hankel matrix in \eqref{eq_Hankelinout}, and let $r:=\text{rank}(H)$. Consider unitary matrices $U=\left[\begin{array}{cc} U_1 & U_2  \end{array}\right]$ and $V=\left[\begin{array}{cc} V_1 & V_2  \end{array}\right]$  that ``compress'' the matrix $H$ as follows:
\begin{equation} \label{eq.compressH1} U^* H V =\left[\begin{array}{cc} 0 &  0 \\ 0 & U_2^* H V_2 \end{array}\right]=: \left[\begin{array}{cc} 0 &  0 \\ 0 & \widehat H \end{array}\right],
\end{equation}
where $\widehat H$ is of dimension $r\times r$ and invertible. Now, we set:
$$\widetilde{L}_{\la}(z):=\left[\begin{array}{cc|c} U^* & &  \\ & V^* &  \\ \hline  &  & I_m \end{array}\right] \widehat L_{\la}(z) \left[\begin{array}{cc|c} U & &  \\ & V &  \\ \hline  &  & I_m \end{array}\right],$$
where $\widehat{L}_{\la}(z)$ is the linear system matrix in Theorem \ref{th_linear_palind_Hankel}. Then, $\widetilde{L}_{\la}(z)$ is a ``compressed'' pencil of the form
\begin{equation}\label{compressed}
\widetilde{L}_{\la}(z)=:\left[\begin{array}{cc;{2pt/2pt}cc|c} 0 & 0 & 0 &  0 &  0 \\
0 & 0 & 0 &  A_c(z) &  F_c  \\ \hdashline[2pt/2pt]
0 & 0 & 0 &  0 &  0 \\
0 & { \rev A_c^* (z)} & 0 & 0& (1+z) G_c \\ \hline 0 & z F_c^* & 0& (1+z) G_c^*  & (1+z) R_0  \end{array}\right],
 \end{equation}
{ where $A_c (z) \in \mathbb{C}[z]^{r\times r}$ is a regular pencil,} and
\begin{equation} \label{2compressed}  L_c(z):= \left[\begin{array}{cc|c}
 0 &   A_c(z) &  F_c  \\
 { \rev A_c^* (z)}  & 0& (1+z) G_c \\ \hline  z F_c^* & (1+z) G_c^*  & (1+z) R_0  \end{array}\right]
\end{equation}
is a strongly minimal $*$-palindromic linearization for $(1+z)R(z)$.
\end{theorem}
\begin{proof} We define the following pencils { based on submatrices of \eqref{pal_pencil_anydegree}}:
\begin{equation*}
\begin{bmatrix}
X_1(z) \\
X_2(z)
\end{bmatrix} :=	\left[\begin{array}{c} K_{\la}(z) \\ \hline \phantom{\Big|} (z+1)(e_d^T\otimes I_m)
	\end{array}\right] \,\,\text{and}\,\, \begin{bmatrix}
	Y_1(z) & Y_{2}(z)
	\end{bmatrix}:= \left[\begin{array}{c|c} K_{\la}^T(z) &  e_d\otimes I_m
	\end{array}\right].
\end{equation*}
Then, we have that
	\begin{equation}\label{eq:rel1}
	\begin{bmatrix}
	X_1(z) \\
	X_2(z)
	\end{bmatrix} H= \begin{bmatrix}
	M_{\la}(z) \\
	(1+z) G^*
	\end{bmatrix}\quad  \text{ and } \quad H \begin{bmatrix}
	Y_1(z) & Y_{2}(z)
	\end{bmatrix}=  \begin{bmatrix}
	M_{\la}(z) & F
	\end{bmatrix}.
	\end{equation}
Since $\begin{bmatrix}
X_1(z) \\
X_2(z)
\end{bmatrix}$ and $\begin{bmatrix}
Y_1(z) & Y_{2}(z)
\end{bmatrix}$ have full column rank and full row rank, respectively, for all $z_0\in\mathbb{C}$ and also at $\infty$, i.e., the matrix coefficients of $z$ of these two pencils have full rank, then
 $\begin{bmatrix}
	M_{\la}(z) \\
	(1+z) G^*
	\end{bmatrix} $ and $ \begin{bmatrix}
	M_{\la}(z) & F
	\end{bmatrix}$ have rank $r$ for all $z_0\in\mathbb{C}$ and at $\infty$. Moreover, the right null space of $\begin{bmatrix}
	M_{\la}(z) \\
	(1+z) G^*
	\end{bmatrix} $  is spanned by the columns of $V_1$ and the left null space of $\begin{bmatrix}
	M_{\la}(z) & F
	\end{bmatrix}$ is spanned by the rows of $U_1^*$. Analogously, we set:
\begin{equation*}
\begin{bmatrix}
\widetilde X_1(z) \\
\widetilde X_2(z)
\end{bmatrix} :=	\left[\begin{array}{c} (\rev K_{\la}^*(z))^T \\ \hline \phantom{\Big|}  z(e_d^T\otimes I_m)
\end{array}\right]
\end{equation*}
and
\begin{equation*}	
\begin{bmatrix}
\widetilde Y_1(z) & \widetilde Y_{2}(z)
\end{bmatrix}:= \left[\begin{array}{c|c} \rev K_{\la}^*(z) &  (1+ z)(e_d\otimes I_m)
\end{array}\right],
\end{equation*}
and we obtain that

\begin{equation} \label{eq:rel2}
\begin{bmatrix}
\widetilde X_1(z) \\
\widetilde X_2(z)
\end{bmatrix} H^*= \begin{bmatrix}
\rev M_{\la}^*(z) \\
z F^*
\end{bmatrix}\, \text{ and } \, H^* \begin{bmatrix}
\widetilde Y_1(z) & \widetilde Y_{2}(z)
\end{bmatrix}=  \begin{bmatrix}
\rev M_{\la}^*(z) & (1+z)G
\end{bmatrix}.
\end{equation}
{ Observe that $\begin{bmatrix}
\widetilde X_1(z) \\
\widetilde X_2(z)
\end{bmatrix}$ and $\begin{bmatrix}
\widetilde Y_1(z) & \widetilde Y_{2}(z)
\end{bmatrix}$ have also full column rank and full row rank, respectively, for all $z_0\in\mathbb{C}$ and at $\infty$. Then
$\begin{bmatrix}
\rev M_\lambda^*(z) \\
z F^*
\end{bmatrix}$ and $\begin{bmatrix}
\rev M_\lambda^*(z)  & (1+z) G
\end{bmatrix}$ have rank $r$ for all $z_0\in\mathbb{C}$ and at $\infty$. Moreover, the right null space of  $\begin{bmatrix}
\rev M_\lambda^*(z) \\
z F^*
\end{bmatrix}$ is spanned by the columns of $U_1$ and the left null space of
$\begin{bmatrix}
\rev M_\lambda^*(z)  & (1+z) G
\end{bmatrix}$ is spanned by the rows of $V_1^*$.

Then, \eqref{eq:rel1} and \eqref{eq:rel2} and the observations about the null spaces and the ranks imply the compressed form \eqref{compressed}, and that the first and second block columns and that the first and second block rows of the $3\times 3$ partitioned pencil $L_c(z)$ in \eqref{2compressed} have each full rank $r$ for all $z_0 \in \mathbb{C}$ and at $\infty$. { This implies that $L_c(z)$ is a strongly minimal linear system matrix, provided that $ A(z):=\left[\begin{array}{cc}
0 &   A_c(z)  \\
\rev A_c^* (z)  & 0 \end{array}\right]$ is regular and that $A_c(z)$ and $\rev A_c^* (z)$ have no eigenvalues in common. Indeed, we now} prove that $A_c(z)$ is regular, which implies that $\rev A_c^* (z)$ and, thus, $A(z)$ are regular. We will also prove that $A_c(z)$ and $\rev A_c^* (z)$ do not have eigenvalues in common (finite nor infinite). Observe that
$$
\begin{bmatrix}
0 & 0 \\
0 & A_c (z)
\end{bmatrix} = U^* K_\lambda (z) H V.
$$
Since $K_\lambda (z)$ is invertible for all $z_0 \in \mathbb{C}$, $z_0 \ne \lambda$, and
at $\infty$, we get that $A_c (z) \in \mathbb{C}[z]^{r \times r}$ is invertible for all $z_0 \in \mathbb{C}$, $z_0 \ne \lambda$, and at $\infty$. Therefore $A_c (z)$ is regular and has at most one eigenvalue equal to $\lambda$. The equality
$$
\begin{bmatrix}
0 & 0 \\
0 & \rev A_c^* (z)
\end{bmatrix} = V^* H^* \rev K_\lambda^* (z) U
$$
proves analogously that $\rev A_c^* (z)$ is regular and has at most one eigenvalue equal to $1/\overline{\lambda}$. Since $\lambda$ is inside the unit circle $\la \ne 1/\overline{\lambda}$ and $A_c(z)$ and $\rev A_c^* (z)$ do not have eigenvalues in common.

Since $L_c (z)$ is clearly $*$-palindromic, it only remains to prove that
the transfer function of $L_c(z)$ is $(1+z)R(z)$. For that, we} define the rational matrix
$$N(z):=\left[\begin{array}{ccc;{2pt/2pt}ccc|c}
\dfrac{(1+z)z^{{d}-1}}{(1-z\overline{\la})^{d}}I_m & \cdots & \dfrac{1+z}{1-z\overline{\la}}I_m & \dfrac{1}{(z-\la)^{d}} I_m & \cdots & \dfrac{1}{z-\la} I_m & I_m
\end{array}\right],$$
which satisfies that
$$ \left[\begin{array}{cc;{2pt/2pt}cc|c} 0 & 0 & 0 &  0 &  0 \\
0 & 0 & 0 &  A_c(z) &  F_c  \\ \hdashline[2pt/2pt]
0 & 0 & 0 &  0 &  0 \\
0 & { \rev A_c^* (z)} & 0 & 0& (1+z) G_c  \end{array}\right]\begin{bmatrix}
U^{*} & 0 & 0 \\
0 & V^{*} & 0\\
0 & 0 & I_m
\end{bmatrix}N(z)^{T}=0.$$
In particular,
$$ \left[\begin{array}{c;{2pt/2pt}c|c}
  0 &   A_c(z) &  F_c  \\ \hdashline[2pt/2pt]
 { \rev A_c^* (z)} & 0& (1+z) G_c  \end{array}\right]
M(z)^{T}=0,$$
where
$  M(z)^{T}:= \begin{bmatrix}
[0 \quad I_r]U^{*} & 0 & 0 \\
0 & [0 \quad I_r]V^{*} & 0\\
0 & 0 & I_m
\end{bmatrix}N(z)^{T}.$ Then, we have by \cite[Lemma 3.1]{dopquinvan2022} that the transfer function of the system matrix $L_c(z)$ is
\begin{align*} \begin{bmatrix}
 z F_c^* & (1+z) G_c^*  & (1+z) R_0
\end{bmatrix} M(z)^{T}
&  = \begin{bmatrix}
z F^*& (1+z) G^* & (1+z) R_0
\end{bmatrix}N(z)^{T} \\ & = (1+z)R(z).\end{align*}
\end{proof}

\begin{remark}[Multiple poles] \rm If  the para-Hermitian rational matrix $R(z)$ has more than one pole inside $S^1$, i.e., $\Omega_{in}=\{\la_1,\ldots,\la_p\}$ and
\begin{equation}\label{eq: rat_mult_pole}
    R(z)= \displaystyle\sum_{i=1}^p R_{\la_i}(z) + R_0+ \displaystyle\sum_{i=1}^p R_{\la_i}^*(1/z) ,
\end{equation}
{ where $R_{\la_i}(z)$ and $R_{\la_i}^*(1/z)$ are as in \eqref{eq:sproper_partial} and \eqref{eq:out_partial},}
we can construct a { strongly minimal} $*$-palindromic linearization from each pole $\la_i$ and combine them, to obtain a strongly minimal $*$-palindromic linearization of $(1+z)R(z)$. More precisely, if $L_{\la_i}(z)$ is a strongly minimal $*$-palindromic linearization of $(1+z)[R_{\la_i}(z)+R_{\la_i}^*(1/z)]$, for $i=1,\ldots,p,$ of the form:
\begin{equation}  L_{\la_i}(z):= \left[\begin{array}{cc|c}
 0 &   P_i(z) &  F_i  \\
 \rev P_i^*(z)  & 0& (1+z) G_i \\ \hline  z F_i^* & (1+z) G_i^*  & 0  \end{array}\right],
\end{equation}
constructed by using Theorem \ref{th:$*$-palindromic_onepole}, Theorem \ref{th_linear_palind_Hankel} or Theorem \ref{th:$*$-palindromic}, then
$$L(z)=\left[\begin{array}{cc;{2pt/2pt}c;{2pt/2pt}cc|c} 0 & P_1(z) &   & &   &  F_1 \\
\rev P_1^*(z) & 0 &  &  &   &  (1+z) { G_1}  \\ \hdashline[2pt/2pt]
 &  & \ddots &  &     &  \vdots \\
\hdashline[2pt/2pt]
 &  &   & 0 &  P_p(z) &  F_p \\
 &    & & \rev P_p^*(z) &  0 &  (1+z) G_p  \\\hline
z F_1^* & (1+z) G_1^* & \cdots & z F_p^* & (1+z) G_p^* & (1+z) R_0
\end{array}\right]$$
is a strongly minimal $*$-palindromic linearization of $(1+z)R(z)$.

\end{remark}

 \begin{remark}[Polar sections at $0$ and $\infty$] \rm We can consider an alternative decomposition for rational matrices that does not require knowledge of their stable and anti-stable parts, based on the polar sections at $0$ and $\infty$. Namely, if a para-Hermitian rational matrix $R(z)\in { \mathbb{C}(z)^{m\times m}}$ is not proper then it must have a pole at $z=\infty$ and a pole at $z=0$, with the same partial multiplicities. Hence, one can decompose $R(z)$ as follows:
\begin{equation} \label{PFE}  R(z) = R_p(z) + R_{0}(z) + R_{\infty}(z), \quad R_{0}(z) = \sum_{i=1}^d R_{-i} z^{-i}, \quad
R_{\infty}(z) = \sum_{i=1}^d R_i z^i,
\end{equation}
where $R_p(z)$ is a proper para-Hermitian rational matrix, and $R_{0}(z)$ and $R_{\infty}(z)$ are the polar sections at $z=0$ and $z=\infty$, respectively, in the Laurent series of $R(z)$ around infinity. The fact that their highest degree coefficient $d$ is the same, follows from the para-Hermitian property. Moreover, this also implies that
$  R_{-i} = R_i^*, \;\; i=1,\ldots,d.$

{ If we} also assume that the proper rational matrix $R_p(z)$ in this decomposition \eqref{PFE} is a constant matrix $R_0$, { then} $R(z)$ is a para-Hermitian rational matrix of the form
\begin{equation}\label{eq_rat_degree}
R(z)=R_dz^d+\cdots+R_1 z+ R_0 + R_{-1}\dfrac{1}{z}+\cdots+ R_{-d}\dfrac{1}{z^{d}},
\end{equation}
{ with} $R_{i}^*=R_{-i}$, for $i=1,\ldots,d$, and $R_{0}^*=R_0$.  Then, if $R_{-d}$ is invertible, we can apply Theorems \ref{th:$*$-palindromic_onepole} or \ref{th_linear_palind_Hankel} { with $\lambda = 0$} to construct a strongly minimal $*$-palindromic linearization of $(1+z)R(z)$. If $R_{-d}$ is not invertible, we can apply Theorem \ref{th:$*$-palindromic} { with $\lambda = 0$}.
 \end{remark}

\begin{example} \rm We consider here the representation in \eqref{eq_rat_degree} in the simplest case $d=1$. That is, $R(z)$ is a $m\times m$ para-Hermitian rational matrix of the form
$$R(z)=R_1 z+ R_0 + R_{-1}\dfrac{1}{z},$$
meaning that $R_{1}^*=R_{-1}$ and $R_{0}^*=R_0$. Let us first assume that $R_{-1}$ (and hence also $R_1$) is invertible. Then, the following linear system matrix $L(z)$ is a strongly minimal $*$-palindromic linearization for $(1+z)R(z)$:
\begin{align*}
 L(z)= \left[\begin{array}{cc|c} 0 & -z I_m & R_{-1} \\ - I_m & 0 & I_m (1+z) \\ \hline z R_1 & I_m (1+z) & R_0(1+z) \end{array}\right].
\end{align*}
{ Observe that $L(z)$ above can be obtained from the pencil \eqref{pal_pencil_anydegree} taking $d=1$ and $\lambda=0$.}
 If $R_{1}$ is not invertible, i.e., $\rank{R_{1}}=r<m$. Then, we can write $R_{1}=L_{1}U_{1}^*,$ with $L_{1},U_{1}\in\C^{m\times r}$ { and $\rank L_1 = \rank U_1 = r$}. And, since $R_{1}^*=R_{-1}$, we have that $R_{-1}=U_1 L_1^*$. We can thus construct the following trimmed pencil
\begin{align*}
\widetilde L(z)= \left[\begin{array}{cc|c} 0 & -z I_r & L_1^* \\ - I_r & 0 & U_1^* (1+z) \\ \hline z L_1  & U_1 (1+z) & R_0(1+z) \end{array}\right] ,
\end{align*}
that is a strongly minimal $*$-palindromic linearization for $(1+z)R(z)$.
\end{example}

\section{Conclusions and future work} \label{sec.conclusions}
Given a para-Hermitian (resp., para-skew-Hermitian) rational matrix $R(z)$, in this paper we show how to construct a strongly minimal $*$-palindromic (resp., $*$-anti-palindromic) linearization for $H(z):=(1+z)R(z),$ whose eigenvalues preserve the symmetries of the zeros and poles of $R(z)$, { and its minimal indices preserve the equality of the left and right minimal indices in the singular case}. In some cases, the proposed techniques require some computations to be performed to construct the linearization; these computations can be performed by using reliable tools such as the SVD. However, a full analysis of the possible stability of the method is left as an open problem.
{ To obtain our main results, we develop several other results on para-Hermitian and para-skew-Hermitian matrices that are also interesting by themselves, as the structured properties of several Möbius transforms acting on these classes of matrices, the properties of additive decompositions into stable and anti-stable parts, and ways to parameterize explicitly para-Hermitian and para-skew-Hermitian matrices.}


\begin{thebibliography}{01}

\bibitem{AhmadMehrmann} Sk. S.~Ahmad, V.~Mehrmann,
\textit{Backward errors for eigenvalues and eigenvectors of Hermitian, skew-Hermitian, H-even and H-odd matrix polynomials}, Linear and Multilinear Algebra, 61(9) (2013), 1244–1266
	
\bibitem{AmMaZa15} A.~Amparan, S.~Marcaida, I.~Zaballa,
\textit{Finite and infinite structures of rational matrices: a local approach}, Electron. J. Linear Algebra, 30 (2015), 196--226.

\bibitem{Antoulas} A.~C. Antoulas, \textit{Approximation of Large-Scale Dynamical Systems}, SIAM, Philadelphia, 2005.

\bibitem{BN} G.~Barbarino, V.~Noferini,
\newblock \textit{On the Rellich eigendecomposition of para-Hermitian matrices and the sign characteristics of palindromic matrix polynomials},
\newblock  Linear Algebra Appl., 672 (2023), 1--27.



\bibitem{DMQVD} F. M.~Dopico, S.~Marcaida, M. C.~Quintana, P.~Van Dooren,
\newblock \textit{Local linearizations of rational matrices with application to rational approximations of nonlinear eigenvalue problems},
\newblock  Linear Algebra Appl., 604 (2020), 441--475.

\bibitem{DNZ} F. M. Dopico, V. Noferini, I. Zaballa,
\newblock \textit{Rosenbrock's theorem on system matrices over elementary divisor domains},
\newblock submitted, available at {\tt https://arxiv.org/abs/2406.18218}.


\bibitem{DQV} F. M.~Dopico, M. C.~Quintana, P.~Van Dooren,
\newblock \textit{Linear system matrices of rational transfer functions},
\newblock``Realization and Model Reduction of Dynamical
	Systems. A Festschrift to honor the 70th birthday of Thanos Antoulas'', pp. 95-113, Springer-Verlag (2022).

\bibitem{dopquinvan2022} F. M. Dopico, M. C. Quintana, P.~Van Dooren, \textit{Strongly minimal self-conjugate linearizations for polynomial and rational matrices}, SIAM J. Matrix Anal. Appl., 43(3) (2022), 1354--1381.

\bibitem{For75} G. D.~Forney,
\newblock \textit{Minimal bases of rational vector spaces, with applications to multivariable linear systems},
\newblock SIAM J. Control, 13 (1975), 493--520.

\bibitem{Geninetal} Y. Genin, Y. Hachez, Y. Nesterov, R. Stefan, P. Van Dooren, S. Xu,
\newblock \textit{Positivity and linear matrix inequalities},
\newblock Eur. J. Control, 8(3) (2002), 275--298.

\bibitem{real}C. Heij, A. Ran, F. van Schagen, \textit{Introduction to Mathematical Systems Theory: Linear Systems, Identification and Control}, Birkh\"{a}user Verlag, Basel, 2007.

\bibitem{Kai80} T.~Kailath,
\newblock \textit{Linear Systems},
\newblock Prentice Hall, Englewood Cliffs, NJ, 1980.


\bibitem{KressnerQRpal} D. Kressner, C. Schr\"{o}der, D. Watkins,
\newblock \textit{Implicit QR algorithms for palindromic and even eigenvalue problems},
\newblock Numer. Algor., 51 (2009), 209–238.



\bibitem{GoodVibra} D. S. Mackey, N. Mackey, C. Mehl, V. Mehrmann, \textit{Structured polynomial eigenvalue problems: Good vibrations from good linearizations}, SIAM J. Matrix Anal. Appl., 28 (2006), 1029--1051.

\bibitem{antitriangular} D. S. Mackey, N. Mackey, C. Mehl, V. Mehrmann, \textit{Numerical methods for palindromic eigenvalue problems: Computing the anti‐triangular Schur form}, Numer. Linear Algebra Appl., 16(1)
(2009), 63-86.


\bibitem{M4Mob} D. S. Mackey, N. Mackey, C. Mehl, V. Mehrmann, \textit{M\"{o}bius transformations of matrix polynomials}, Linear Algebra Appl., 470 (2015), 120–184.


\bibitem{Nof12} V. Noferini, \textit{The behaviour of the complete eigenstructure of a polynomial matrix under a generic rational transformation}, Electron. J. Linear Algebra, 23 (2012), 607–624.

\bibitem{NV23} V. Noferini, P. Van Dooren, \textit{Root vectors of polynomial and rational matrices: theory and computation}, Linear Algebra Appl.,  656 (2023), 510--540.

\bibitem{Rosen} H. H. Rosenbrock,
\textit{State-space and Multivariable Theory,} Thomas Nelson and Sons, London, 1970.


\bibitem{vandooren-laurent-1979} P. Van Dooren, P. Dewilde, J. Vandewalle, \textit{On the determination of the Smith-McMillan form of a rational matrix from its Laurent expansion}, IEEE Trans. Circuit Syst., 26(3) (1979), 180--189.



\bibitem{WPP} S.~Weiss, J.~Pestana, I.~K. Proudler,
\newblock \textit{On the existence and uniqueness of the eigenvalue decomposition of a parahermitian matrix},
\newblock  IEEE Trans. Signal Process., 66(10) (2018),  2659--2672.



\end{thebibliography}
\end{document}